\documentclass{article}
\usepackage{amssymb,amsmath,amsfonts,enumerate,url,graphicx}

\newtheorem{theorem}{Theorem}[section]

\newtheorem{definition}[theorem]{Definition}
\newtheorem{example}[theorem]{Example}

\begin{document}

\title{Envelopes of Circles Centered on a Kiss Curve}

\author{Thierry Dana-Picard and Daniel Tsirkin}

\maketitle

\begin{center}
Department  of  Mathematics,  Jerusalem  College  of
Technology\\  
Jerusalem 9116011, Israel\\e-mail: ndp@jct.ac.il, academictsirkin@gamil.com
\end{center}

\textbf{Abstract:}
Envelopes of parameterized families of plane curves is an important topic, both for the mathematics involved and for its applications. Nowadays, it is generally studied in a technology-rich environment, and automated methods are developed and implemented in software. The exploration involves a dialog between a Dynamic Geometry System (used mostly for interactive exploration and conjectures) and a Computer Algebra System (for algebraic manipulations). We study envelopes of families of circles centered on the so-called kiss curve and offsets of this curve, observing the differences between constructs. Both parametric presentations and implicit equations are used, switching from parametric to polynomial representation being based on packages for Gr\"obner bases and Elimination. Singular points, both cusps and points of self-intersection (crunodes), are analyzed.


\maketitle

\section{Introduction}
\label{intro}
The study of envelopes of 1-parameter families of plane curves disappeared from the undergraduate syllabus, quite early in the 2nd half of the 20th century. Thom complained about that in \cite{thom}, explaining that maybe the reason is that this theory has too few theorems and too many special cases. Nevertheless, their importance both in mathematics and applied fields has always been recognized (e.g., see \cite{bickel}) and they appeared as chapters in books such as \cite{berger,bruce and giblin}, monographs \cite{arnold} and papers \cite{revival} have continuously been published. Some reveal new aspects of classical topics, such as plane curves \cite{strings-envelopes - ACTM 2023,cassini ovals} or applications such as safety zones in industrial plants and entertainment parks \cite{safety}. 

Offsets of algebraic plane curves, also called parallel curves, have been studied by Leibnitz \cite{leibnitz}; since then they drew continuous interest. Among the difficulties we shall mention that their study involves heavy computations, the degree of offset equation is generally much larger than that of the original curve (called the \emph{progenitor}) and the topology of the offset can be much more complicated than that of the original curve. After the above remarks, Alcazar and Sendra mention in \cite{alcazar-sendra} that if the original curve is the Folium of Descartes (whose equation is $x^3+y^3-3xy=0$), the polynomial defining an offset is of degree 14 and has 114 terms. This question has been addressed by San Segundo and Sendra \cite{san segundo - sendra}.

The envelope of a 1-parameter family of circles of constant radius $r$ centered on a plane curve $\mathcal{C}$  and the offset at distance $r$ of $\mathcal{C}$ (also called a parallel curve) are two different objects. Sometimes they are equal, but sometimes there is a strict inclusion, the smaller being the offset.  

\begin{figure}[htb]
\begin{center}
\includegraphics[width=5cm]{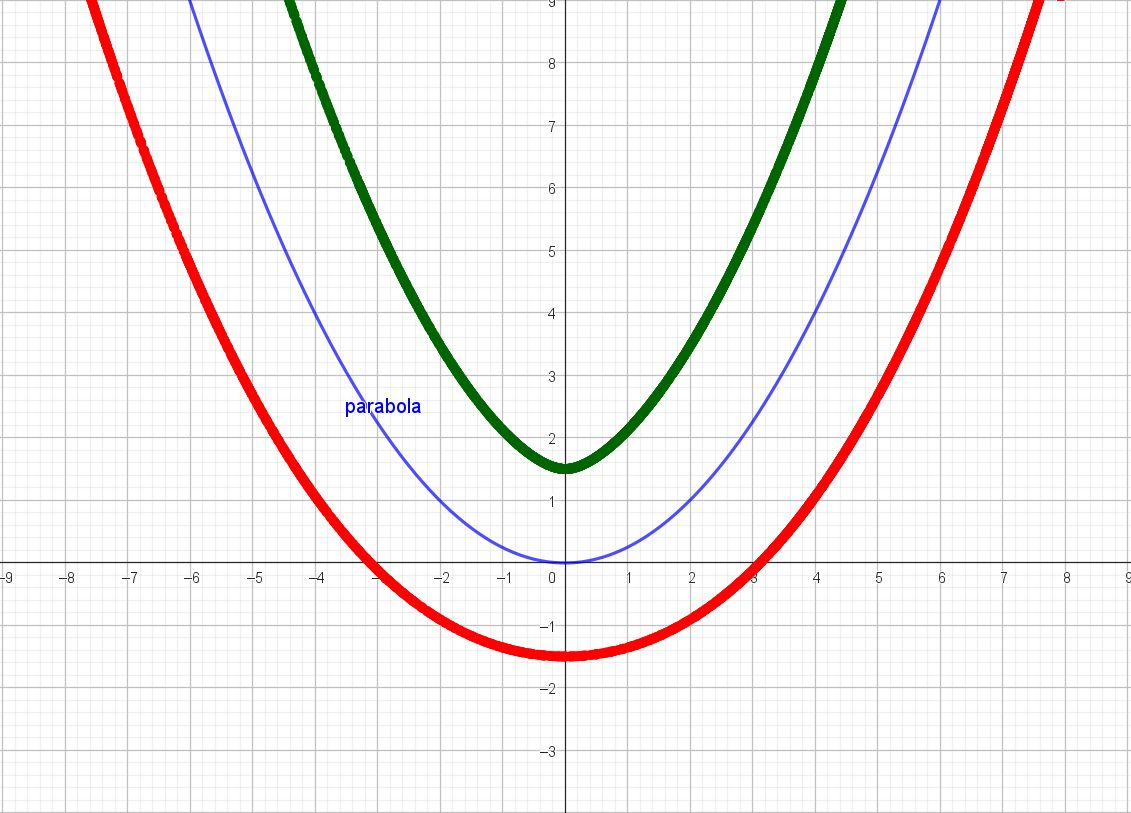}
\qquad
\includegraphics[width=5cm]{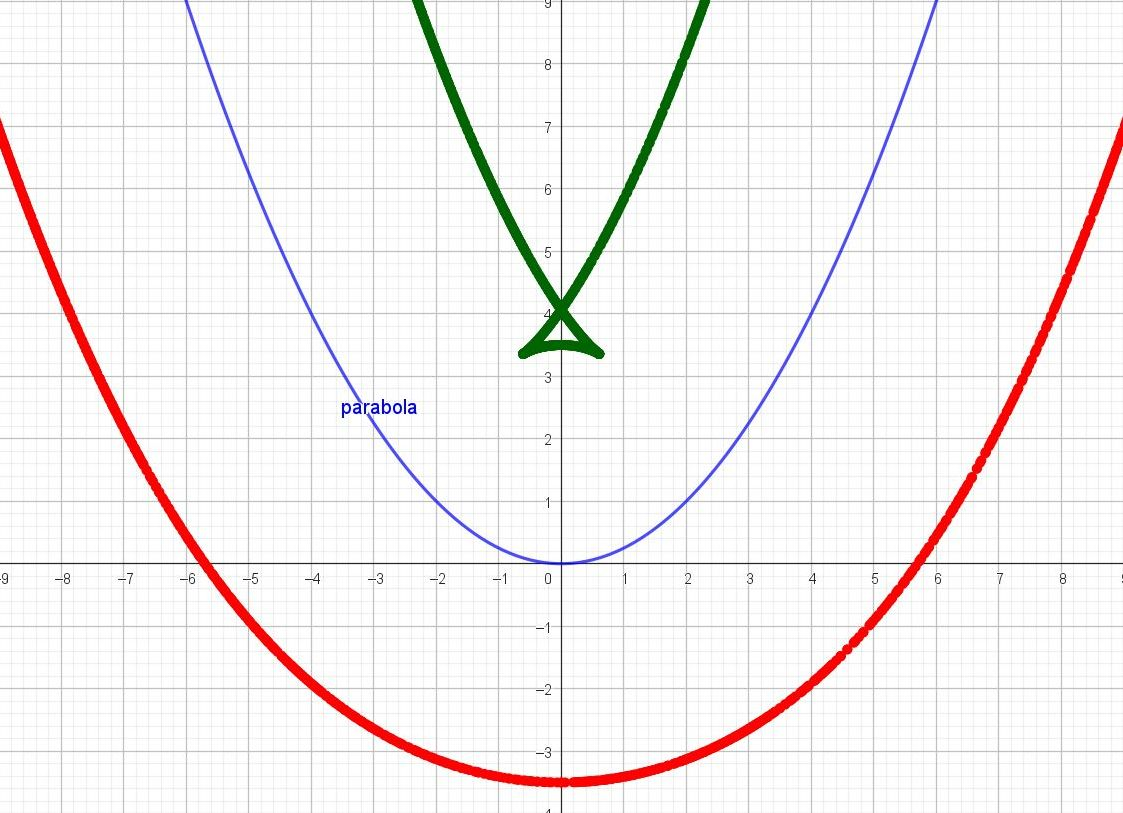}
\caption{Offsets at distance 1.5 and 3.5 of the parabola $\mathcal{C}: 4y-x^2=0$. }
\label{offset parabola}
\end{center}
\end{figure}

As an easier example, Figure \ref{offset parabola} shows offsets of the parabola $\mathcal{C}: 4y=x^2$ for at different distances (a formal definition is given in subsection \ref{subsection offsets}. The two components of the offset at distance 3/2 explain why Leibnitz called this parallel curves, the other example shows why the notion is actually not intuitive. It incites to study the singular points of the offset. A nice 5-vertices star has been constructed in a similar way from a parabola, and analyzed in \cite{solitaria}.

The study of offsets is in fact a study of geometric loci, which is a classical topic in Geometry. The heavy algebraic machinery needed for offsets, and also for envelopes,  may be out of reach of hand-made computations, but also of Computer Algebra Systems. Therefore, works have been devoted to approximate methods for their determination \cite{schutz-juttler,elber et al.}. This is beyond the scope of this paper, and we present situations where algebraic methods work well. For general studies of the topology of offsets, we refer to \cite{alcazar-sendra}, who study various aspects of the offset from the original curve,  and to \cite{seong et al} for their study of self-intersections of offsets. 

In subsections \ref{subsection envelopes of 1-param fam} and \ref{subsection offsets}, we recall the basic definitions of envelopes and offsets. In Section \ref{section networking}, we give a brief survey on the usage of mathematical software, namely Computer Algebra System (CAS) and Dynamic Geometry Software (DGS). They have different but complementary abilities, whence the necessity of  a dialog between them. Then, in Section \ref{section kiss curve}, we study envelopes of families of circles centered on a kiss curve and offsets of this curve. We give detailed computations and compare the two kinds of curves obtained.
Special attention is given to the differences between the curves and to their singular points in certain cases, no exact answer could be obtained, only numerical answers. In our study, we follow an experimental process  in a technology rich environment, and emphasize the segment Exploration-Discovery-Conjecture-Proof. Animation are central for the exploration, and automated methods are almost ubiquitous. 
 
\subsection{Envelope of a 1-parameter family of plane curves}
\label{subsection envelopes of 1-param fam}
Bruce and Giblin \cite{bruce and giblin} give 4 different, non equivalent, definitions of an envelope and describe who is a subvariety of who.  The same definitions appear in \cite{kock,bickel} with different names. It seems that there is still no standard taxonomy. The names of the definitions below are taken from \cite{kock}.  
\begin{definition} [Synthetic-limit definition]
\label{def envelope limit}
Let $\mathcal{C}_k$ be a family of plane curves. The characteristic $\mathcal{M}_k$ is the limit of  $\mathcal{C}_{k+h} \cap \mathcal{C}_k$ when $h \rightarrow 0$. The union of the characteristics is called an envelope of $\mathcal{C}_k$.
\end{definition}

\begin{definition}[Impredicative]
$\mathcal{E}$ is a curve with the property that at each of its points, it is tangent to a unique curve from the family $\mathcal{C}_k$. (Also, $\mathcal{E}$ should touch each $\mathcal{E}$.)
\end{definition}

\begin{definition}[Analytic]
\label{def envelope analytic}
Let $\mathcal{C}_k$ be a family of plane curves given by the equation $F(x,y,k)=0$, where $k$ is a real parameter. An envelope of this family, if it exists, is determined by the solutions of the system of equations:
\begin{equation}
\label{eq envelope}
\begin{cases}
F(x,y,k)=0\\
\frac{\partial F}{\partial k} F(x,y,k)=0
\end{cases}.
\end{equation}
\end{definition}

Together with an elementary example, \cite{revival} shows that Definition \ref{def envelope limit} implies Definition \ref{def envelope analytic}. A general proof is to be found in  \cite{bruce and giblin}. Berger \cite{berger} gives only the analytic definition.

\begin{example}
The above mentioned example of \cite{revival} is as follows: a family of lines is given in the plane by the equation $x+ky=k^2$, where $k \in \mathbb{R}$.  This family has an envelope, namely the parabola whose equation is $y^2+4x=0$; see Figure \ref{fig env lines}, which has been obtained with GeoGebra\footnote{A multipurpose software, freely downloadable from \url{http://www.geogebra.com}. Originally, it was a Dynamic Geometry Software, but now it has many other abilities.} Note the slider bar in the upper left corner; the availability of sliders is a core feature of Dynamic Geometry. A GeoGebra applet is available at \url{https://www.geogebra.org/m/rfewh5xh}.

\begin{figure}[htb]
\begin{center}
\includegraphics[width=5cm]{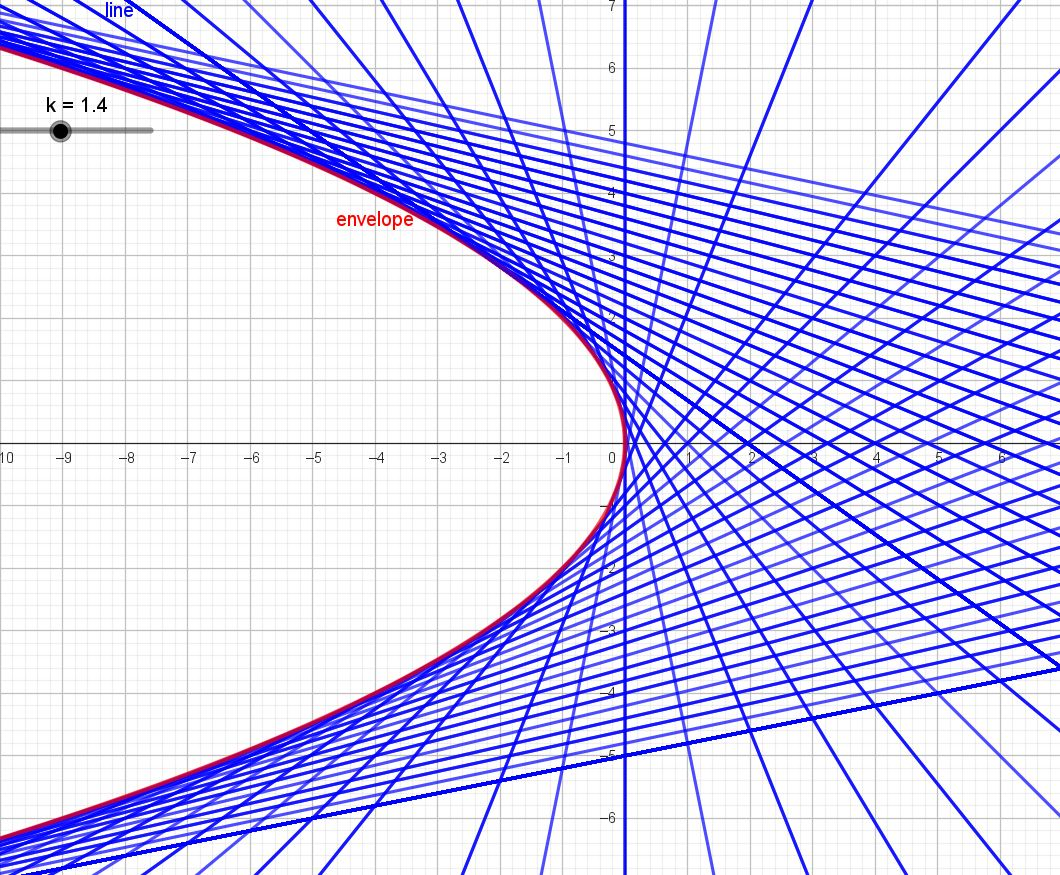}
\caption{Envelope of a family of lines}
\label{fig env lines}
\end{center}
\end{figure}
\end{example}

\begin{example}
\label{example 2}
 Consider the circle centred at the origin with radius 2, and the family of circles centered on it and tangent to the $y-$axis. This family has an envelope, which is a nephroid. Figure \ref{fig nephroid} shows a construction with the automated methods implemented in GeoGebra, namely the command \textbf{Envelope(<path>,<point>)}. This command provides a double output: the plot displayed in Figure \ref{fig nephroid} and a symbolic equation 
\begin{equation}
\label{eq nephroid}
x^7 + 3x^5 y^2 - 12x^5 + 3x^3 y^4 - 24x^3 y^2 - 60x^3 + x y^6 - 12x y^4 + 48x y^2 - 64x = 0.
\end{equation}
Denote by $F(x,y)$ the left hand side in Equation (\ref{eq nephroid}). Actually, we have:
\begin{equation}
\label{eq factorized envelope}
F(x,y)=x \; (x^6 + 3x^4 y^2 - 12x^4 + 3x^2 y^4 - 24x^2 y^2 - 60x^2 + y^6 - 12 y^4 + 48y^2 - 64).
\end{equation}
The right factor can be identified as determining a nephroid, which is a sextic (see \url{https://mathcurve.com/courbes2d.gb/nephroid/nephroid.shtml}). The $x$ factor determines the $y-$axis. It is clear that only a segment of the $y-$axis is tangent to the circles in the family. This segment cannot be determined by algebraic methods only: the software works with computations of Gr\"obner bases, and uses elimination. Therefore it gives a Zariski closure of the envelope. The GeoGebra applet is available at \url{https://www.geogebra.org/m/qjthpat4}.
 \begin{figure}[htb]
\begin{center}
\includegraphics[width=5cm]{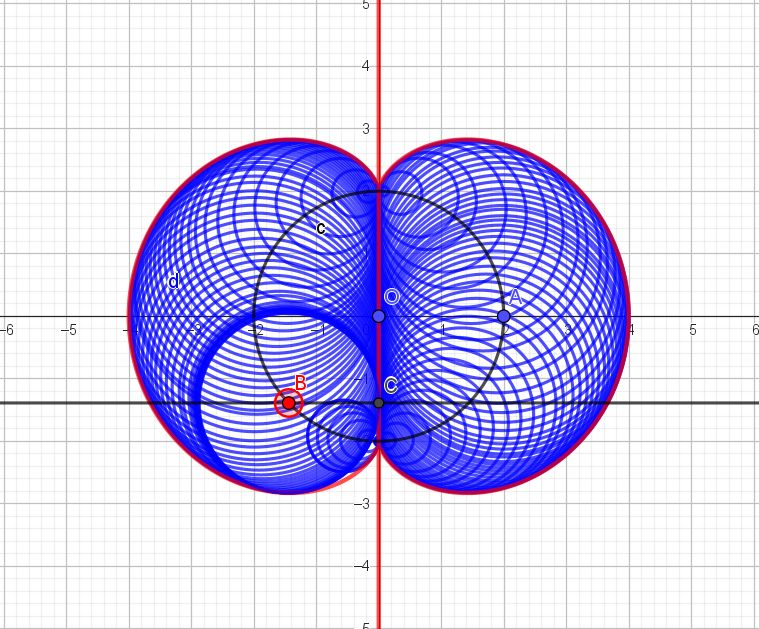}
\caption{A nephroid as an envelope of a family of circles}
\label{fig nephroid}
\end{center}
\end{figure}
\end{example}

\subsection{Offsets}
\label{subsection offsets}
\begin{definition}
Let $\mathcal{C}$ be a plane curve and $d$ be a positive real number. A regular point $A$ of $\mathcal{C}$ is a point where $\mathcal{C}$ is smooth, i.e. tangents vectors and normal vectors are defined. If $A$ is a regular point, we denote by $\overset{\rightarrow}{N_A}$ a unit vector normal to  $\mathcal{C}$ at $A$. Let $B$ be a point such that $\overset{\rightarrow}{AB}=d \cdot \overset{\rightarrow}{N_A}$ and $C$ a point such that $\overset{\rightarrow}{AC}=-d \cdot \overset{\rightarrow}{N_A}$. The union of the geometric loci of $B$ and $C$ when $A$ run over $\mathcal{C}$ is called the offset of  $\mathcal{C}$ at distance $d$.

The curve $\mathcal{C}$ is called the progenitor of the offsets
\end{definition}

Suppose that $\mathcal{C}$ is defined by a polynomial equation $F(x,y)=0$. At a regular point, then we can take as a unit vector in the direction of $\nabla F (P)$. The loci of $B$ (resp. of $C$) is called the external (resp. internal) offset of $\mathcal{C}$ at distance $d$. In \cite{alcazar-sendra}, they are denoted by $\mathcal{O}_{+d}(\mathcal{C})$ (resp.  $\mathcal{O}_{-d}(\mathcal{C})$. The examples that will be observed later show that the actual meaning of the adjectives external and internal is not the intuitive one.

The topology of offsets is very different form that of envelopes. It has been explored for numerous plane curves, for example in  \cite{alcazar-sendra,Vrsek-Lavicka,dialog,safety}. We propose here a new contribution; see subsection \ref{subsection offsets}.

GeoGebra has  an automated command for the determination of geometric loci, with numerous versions. The output is generally a plot in the geometric window, together with an equation in the algebraic window. It is well known that the two windows are fully synchronized.
 
\begin{example}
We consider an ellipse $\mathcal{E}$ whose equation is $x^2/25 + y^2/16=1$. Figure \ref{offsets of an ellipse}, obtained with GeoGebra,  shows offsets at distance 1 and 3 respectively of $\mathcal{E}$.
\begin{figure}[htb]
\begin{center}
($d=1$)\includegraphics[width=5cm]{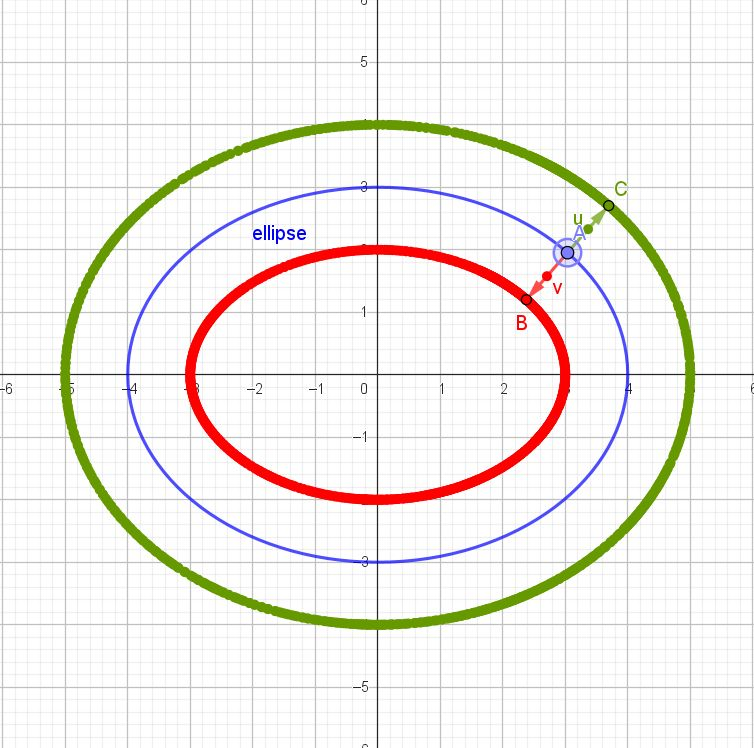}
\qquad
($d=3$)\includegraphics[width=5cm]{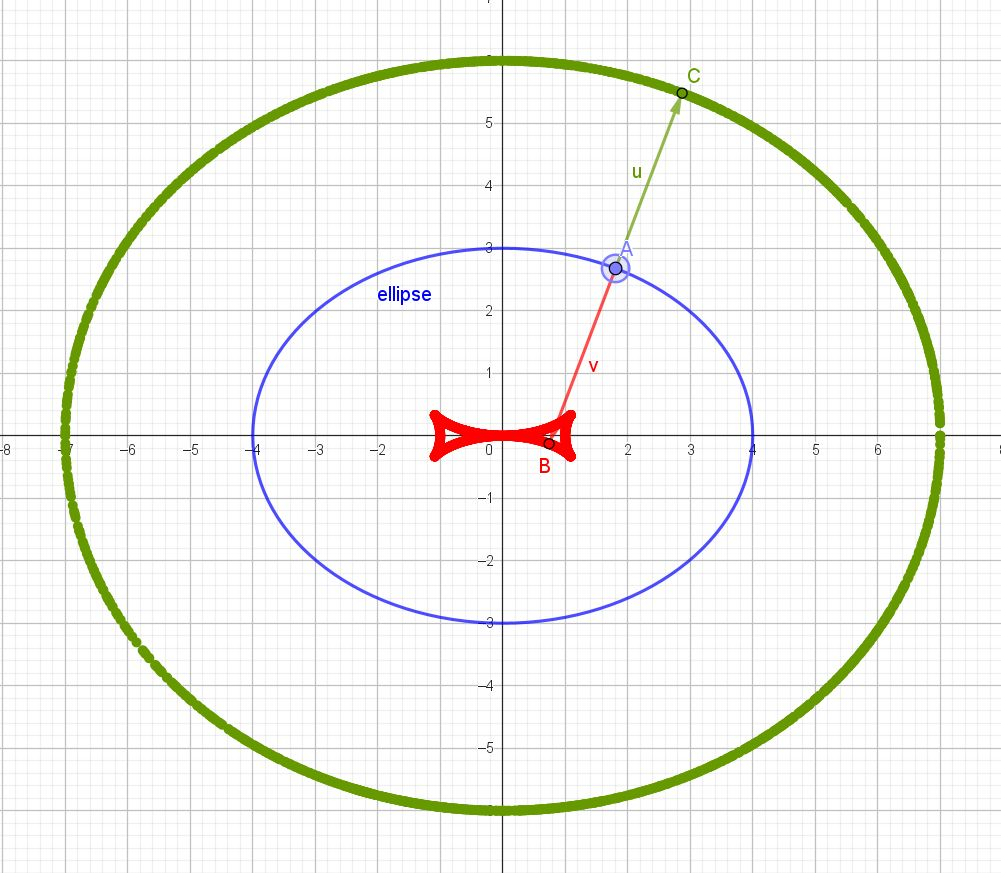}
\caption{Curves parallel to an ellipse}
\label{offsets of an ellipse}
\end{center}
\end{figure}
Note that at distance 1, the offset seems to be smooth, and at distance 3, the offset seems to have 4 cusps. It is easy to understand why they are symmetric about the coordinate axes. Nevertheless, a formal proof is needed for the existence of cusps. This is an example of a general feature of offsets: the fact that the progenitor is smooth does not ensure that an offset will be smooth. On the contrary, the existence of singular points on the progenitor has a strong influence on the topology of an offset.   
\end{example}

Other examples are displayed in \cite{revival}. In \cite{safety}, envelopes of families of circles centred on an astroid are studied. These envelopes are supersets of offsets; their behaviour in a neighborhood of a cusp of the astroid is explored, using zooming offered by the DGS. Such an exploration is important, as the outcome is not intuitive.

\section{Exploration within networking with software}
\label{section networking}
The methods used for the two examples in Section \ref{intro} are different. For the 1st example, computations have been performed using a Computer Algebra System, but they are simple enough to be performed by hand (this has been shown in \cite{revival}). 

For Example \ref{example 2}, work has been performed using the automated methods implemented in GeoGebra. Nevertheless, the algorithms at work are numerical and may produce a not so exact answer. This has been significantly improved in a companion called GeoGebra-Discovery\footnote{Freely downloadable from \url{https://github.com/kovzol/geogebra-discovery}. We suggest to follow the last available version.}. A basic feature of GeoGebra's dynamic geometry is the necessity to make geometric constructs. This is what we did for our nephroid:
\begin{itemize}
\item First, a circle is built, defining its center $O$ and a point $A$ through which the circle has to pass. 
\item A point $B$ is attached to this circle.
\item A line through $B$ and perpendicular to the $y-$ axis is defined. It intersects $y-$axis in a point $C$.
\item A circle $d$ is centered at $B$ and passes through $C$. It is tangent to the $y-$axis.
\item The study explores the envelope of the family of circles $d$ when $B$ runs on the circle $c$. 
\item The \textbf{Move} command is applied to $B$ with \textbf{Trace On} for $d$.
\item The \textbf{Envelope} command is applied, as described above.
\end{itemize}

A more classical way is to solve the suitable system of equations of the form (\ref{eq envelope}); it has been done in \cite{revival}.

Not every situation can be provided by a pure geometric construct. In such a case, GeoGebra's automated command \textbf{Envelope} does not work. GeoGebra can still be used to explore the given situation and to conjecture the existence of an envelope. The exploration uses intensively the dynamic features: sliders, which enable useful animations, and dragging points with the mouse. We must mention that sometimes a slider can be replaced by a segment, the needed parameter being the length of the segment. Nevertheless, this does not allow an animation without dragging.

\section{The Kiss Curve and circles centered on it}
\label{section kiss curve}
The situation that we explore in the present paper is based on the Kiss Curve (a name coined by R. Ferr\'eol in a recent addition to his website Mathcurve \cite{mathcurve}.). This curve is determined by a parametric representation, and, as far as we know, it cannot be constructed by a way acceptable by GeoGebra's command \textbf{Locus}. 

The curve $\mathcal{C}$ is given by a trigonometric parametric representation. This enable to give a symbolic (non polynomial) equation of circles centered on the curve  $\mathcal{C}$ with a given radius. In order to explore the existence of an envelope, we solve the system of equations (\ref{eq envelope}). This yields a trigonometric parametric representation of two components, who are separately and together envelopes of the family. The solution of the system of equations requires heavy algebraic machinery, for which we use the Maple software. 

The exploration of the topology of the envelope, in particular the existence of singular points, requires also heavy machinery. The derivation of a symbolic equation for the envelope needs  a transformation of the equations into a polynomial presentation and the usage of algorithms for algebraic curves (Gr\"obner package and others). This is the topic of next section.

\subsection{The progenitor}
\begin{figure}[htb]
\begin{center}
\includegraphics[width=4cm]{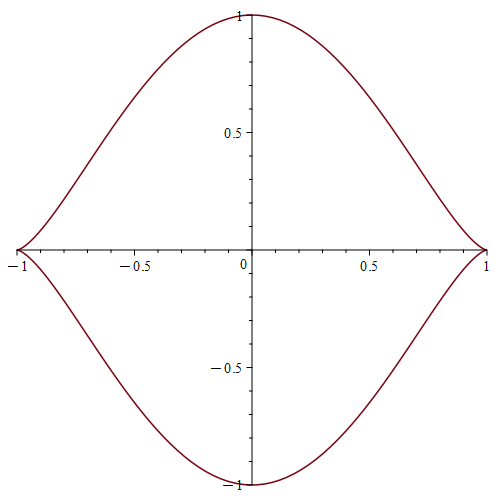}
\caption{A kiss curve}
\label{fig kiss curve}
\end{center}
\end{figure}

The curve $\mathcal{C}$ is given by the following parametric representation:
\begin{equation}\
\label{eq kiss curve}
\begin{cases}
x= \cos u\\
y=  \sin ^3 u
\end{cases}, 
\quad u \in [0, 2\pi ].
\end{equation}
The curve is shown in Figure \ref{fig kiss curve}.
It is easy to prove that the kiss curve is symmetric about the  coordinate axes and has 2 cusps, at coordinates (1,0) and (-1,0). The existence of cusps has a direct influence on the shape of the possible offsets.

\subsection{Envelopes of families of circles}
\begin{figure}[htb]
\begin{center}
($d=1/3$) \includegraphics[width=4cm]{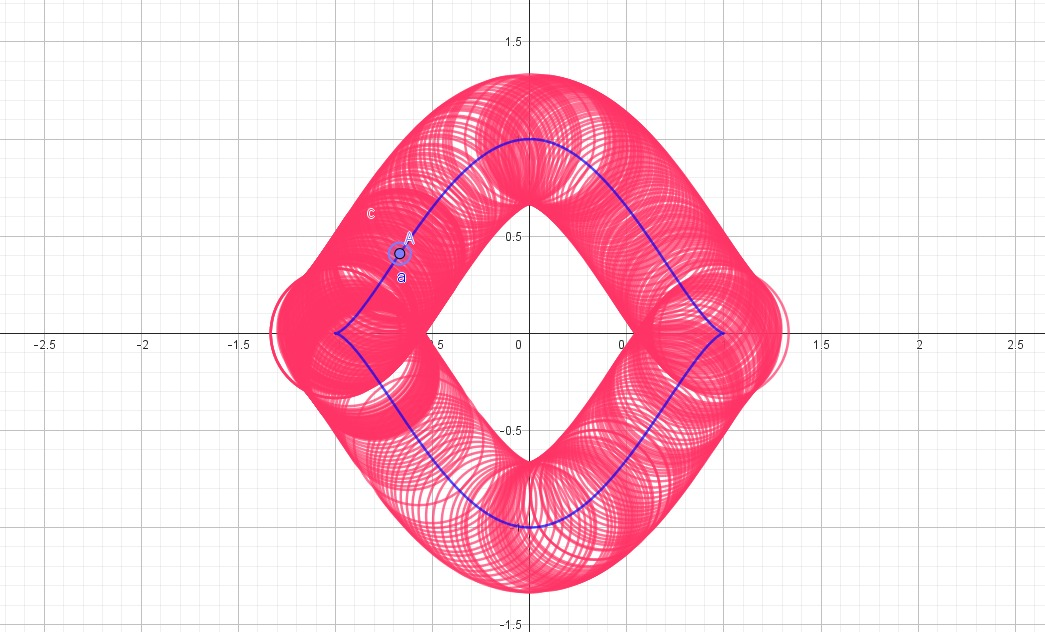}
\qquad 
($d=1/2$)\includegraphics[width=4cm]{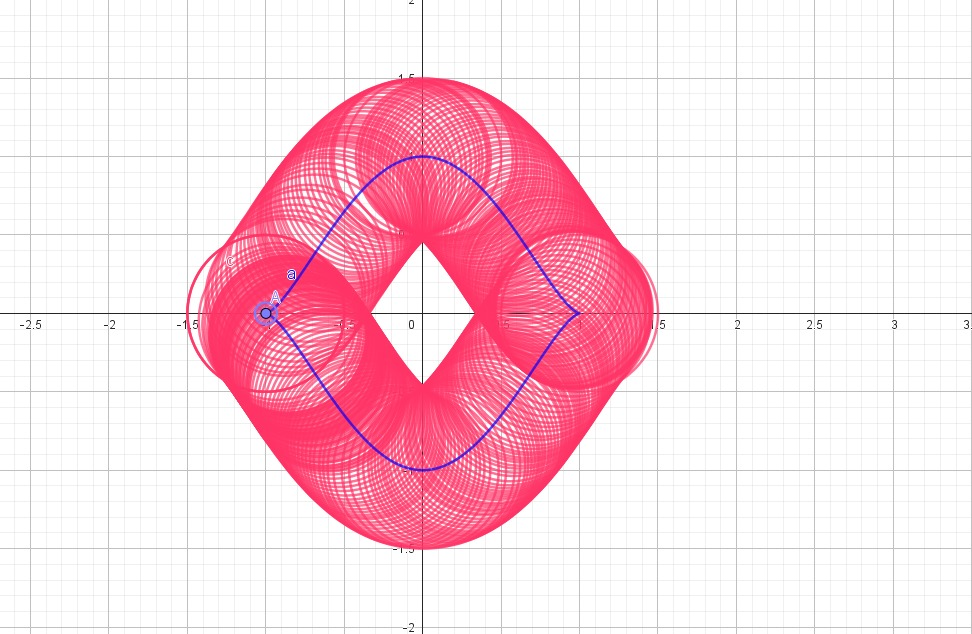}
\qquad 
($d=1$) \includegraphics[width=4cm]{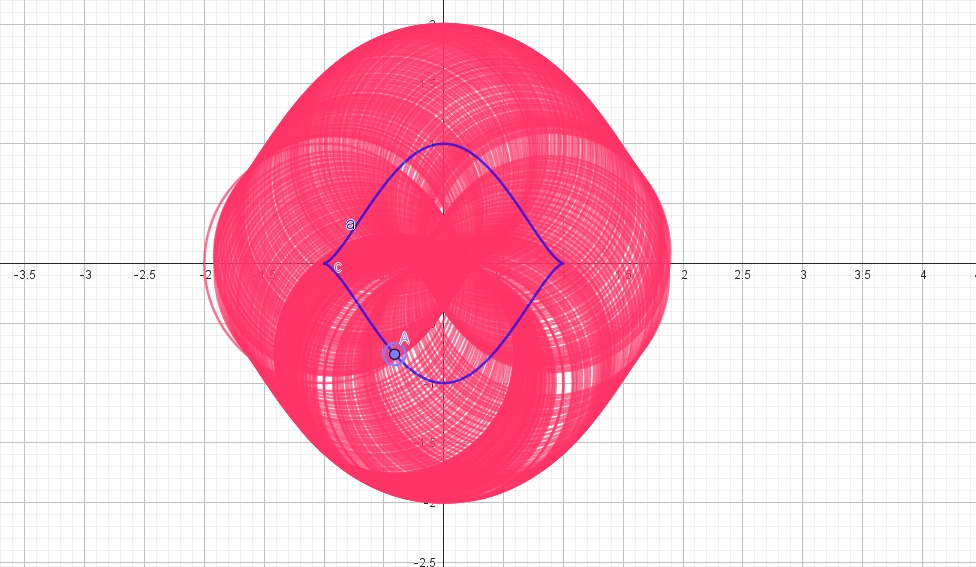}
\qquad
($d=3/2$) \includegraphics[width=4cm]{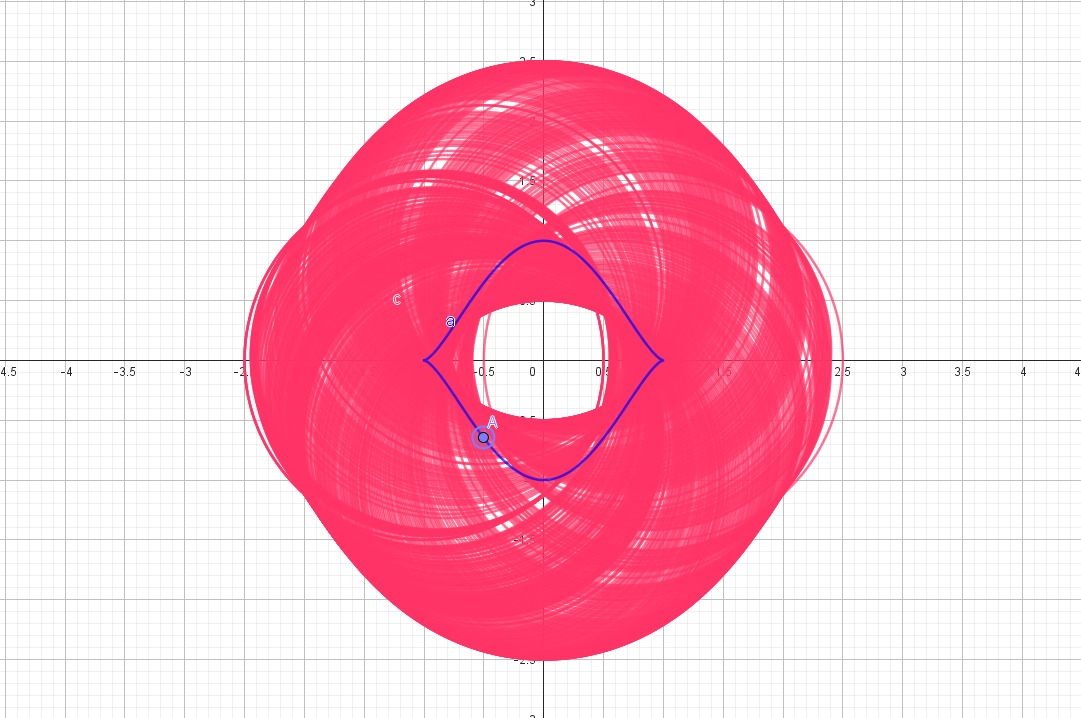}
\caption{Envelopes of families of circles centered on the kiss curve}
\label{fig envelopes}
\end{center}
\end{figure}

We consider a family of circles $\mathcal{F}_r$ centered on $\mathcal{C}$ with radius $k>0$.  Their equation is
\begin{equation}
\label{eq circles}
(x-\cos u)^2-\left( y- \sin^3 u \right)^2-k^2=0.
\end{equation}
Using the DGS, we explore the family and the possible envelope using an animation with \textbf{Trace On}.  Figure \ref{fig envelopes} shows plots for $k=1/3, 1/2, 1$ and $3/2$.  The kiss curve has been plotted with a parametric representation for a point $A$ depending on one parameter, whence the definition of a slider (it appears in the screenshot for $d=3/2$), and then using the command \textbf{Locus($<$point$>$,$<$slider$>$)}. The radius of the circles is encoded by a segment (which appears also in the screenshot for $d=3/2$). The circles have been defined and moved with \textbf{Animation On}. An applet is available at \url{https://www.geogebra.org/m/q45b5rwx}.

\begin{figure}[htb]
\begin{center}
\includegraphics[width=3.5cm]{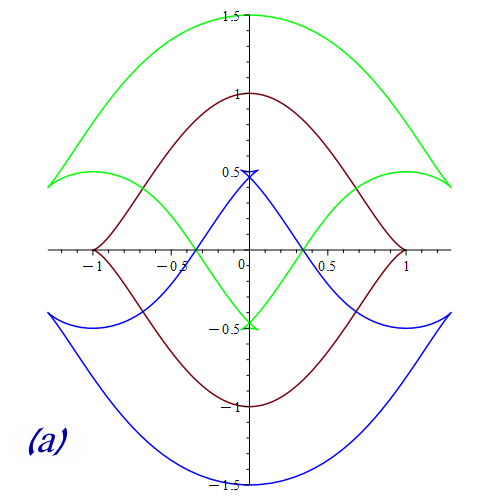}
\quad
\includegraphics[width=3.5cm]{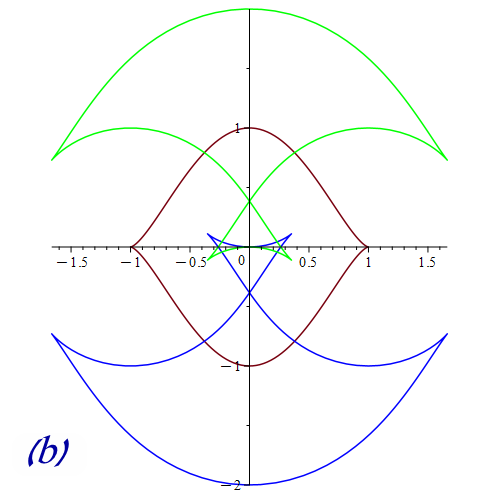}
\quad
\includegraphics[width=3.5cm]{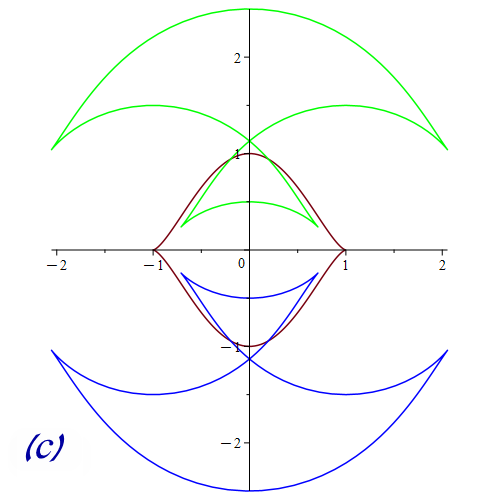}
\caption{The "analytic" envelope of the family of unit circles centered on the kiss curve}
\label{fig envelope of unit circles with Maple}
\end{center}
\end{figure}

Now we wish to determine equations for the possible envelope. 
Denote by $F(x,y,k)$ the left-hand side in Equation (\ref{eq circles}).  According to Definition \ref{def envelope analytic},  the family has an envelope if the following system of equations, has solutions:
\footnotesize
\begin{equation}
\label{eq envelope circles}
\begin{cases}
\left(x-\cos \left(u\right)\right)^{2}+\left(y-\sin \left(u\right)^{3}\right)^{2}-R^{2}=0 \\
2 \sin \left(u\right) \left(\cos \left(u\right) \left(3 \sin \left(u\right)^{4}-3 \sin \left(u\right) y-1\right)+x\right) =0
\end{cases}.
\end{equation}
\normalsize
The output with Maple contains the place holder \emph{RootOf}, which is resolved using the command \textbf{allvalues}. The answer is finally the union of two components, each one given by a parametric representation:
\begin{equation}
\label{1st component envelope}
\begin{cases}
x  = \frac{\cos  u \left(36 \cos^4 u   +3 \sin u  \sqrt{2}\, \sqrt{R^{2} \left(17-9 \cos   4 u \right)}-36 \cos^2 u-4\right)}{36 \cos^4 u-36 \cos^2 u-4},
\\
y =\frac{72 \sin^7   -72 \sin   \left(u \right)^{5}-8 \sin^3 u +2 \sqrt{2}\, \sqrt{R^{2} \left(17-9 \cos 4u)\right)}}{-17+9 \cos   4u}
\end{cases}
\end{equation}
and 
\begin{equation}
\label{2nd component envelope}
\begin{cases}
x = -\frac{3 \left(-12 \cos^4 u  +\sin  u \sqrt{2}\, \sqrt{R^{2} \left(17-9 \cos  4u \right)}+12 \cos^2 u  +\frac{4}{3}\right) \cos u}{36 \cos^4 u  -36 \cos^2 u  -4},
\\
y=\frac{72 \sin^7 u  -72 \sin^5 u  -8 \sin^3 u  -2 \sqrt{2}\, \sqrt{R^{2} \left(17-9 \cos  4u \right)}}{-17+9 \cos  4u}
\end{cases}
\end{equation}
Figure \ref{fig envelope of unit circles with Maple} shows  a snapshot of this envelope, when (a) $R=0.5$, (b) $R= 1$ and (c) $R= 1.5$, obtained with Maple. Figure \ref{fig envelopes} suggests another shape for the outer part for an envelope. This reinforces the fact that an envelope according to Definition \ref{def envelope analytic} may be a strict subvariety of an envelope in the sense of  Definition \ref{def envelope limit}.

\subsection{Looking for singular points of the envelope}
\label{subsection singular points of the envelope}
The experimental work and visual analysis of the output, as shown in Figure \ref{fig envelopes}, suggests that the envelope has singular points, namely several cusps and double points. The determination of cusps follows a well-defined process, but the determination of double points is harder.

This can be explored using the parametric representation of the envelope given in Equations (\ref{1st component envelope}) and (\ref{2nd component envelope}). Because of the symmetries of the situation, it is enough to check the existence of cusps for one of the above components. In what follows, we will work with Equations (\ref{1st component envelope}) and $R=1$.
One of the components of the envelope is given by the equations:
\begin{equation}
\label{eq 1st component}
\begin{cases}
x = \frac{\cos  u \left(36 \cos^4 u  +3 \sin  u \sqrt{2}\, \sqrt{17-9 \cos  4u}-36 \cos^2 u -4\right)}{36 \cos^4 u  -36 \cos^2 u  -4},
\\
y = \frac{72 \sin^7 u  -72 \sin^5 u  -8 \sin^3 u  +2 \sqrt{2}\, \sqrt{17-9 \cos  4u}}{-17+9 \cos  4u}
\end{cases}
\end{equation}

\begin{figure}[htb]
\begin{center}
\includegraphics[width=5cm]{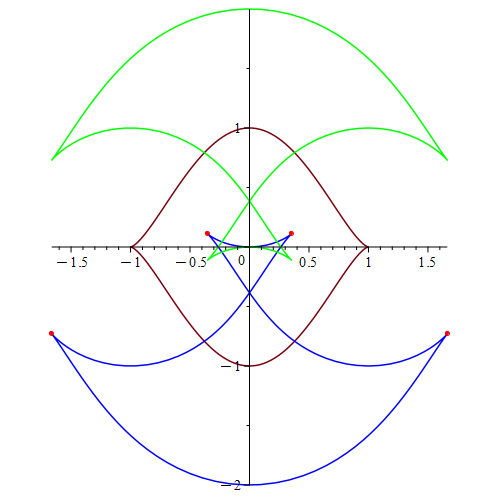}
\caption{The "analytic" envelope with emphasized singular points}
\label{fig envelope with singular points}
\end{center}
\end{figure}

Recall the following definition:
\begin{definition}
Let $\mathcal{C}$ be a plane curve determined by the parametric equations $x=x(u), \; y=y(u)$, where $u$ is a real parameter.
If $x$ and $y$ are differentiable at $u_0$ and if $\frac{\partial x}{\partial u}=\frac{\partial y}{\partial u}=0$, the point $(x(u_0),y=y(u_0))$ is a singular point of $\mathcal{C}$.
\end{definition}
Here we have (we give Maple's as is):
\footnotesize
\begin{equation}
\label{eq 1st derivatives}
\begin{cases}
\frac{dx}{du} = \frac{9 \left(\frac{2 \left(-2 \cos^2 u +1\right) \sqrt{2}}{3}+\sin  u \left(\cos^4 u  -\cos^2 u  -\frac{1}{9}\right) \sqrt{17-9 \cos  4u}\right)}{\sqrt{17-9 \cos  4u}\, \left(9 \cos^4 u  -9 \cos^2 u  -1\right)}\\
\frac{dy}{du} = \frac{27 \cos  u \;  \left(\frac{2 \left(-2 \cos^2 u +1 \right) \sqrt{2}}{3}+\sin  u \left(\cos^4 u  -\cos^2 u  -\frac{1}{9}\right) \sqrt{17-9 \cos  4u}\right) \sin  u}{\sqrt{17-9 \cos  \left(4 u\right)}\, \left(9 \cos^4 u  -9 \cos^2 u  -1\right)}
\end{cases}
\end{equation}
\normalsize
 After simplification, it is enough to look for values of the parameter $u$ for which the numerators vanish simultaneously; we solve the following equations:
 \footnotesize
 \begin{equation}
 \label{eq vanishing numerators}
 \begin{cases}
 \left(-9 \cos^4 u +9 \cos^2 u +1\right) \sin u  \sqrt{17-9 \cos  4u}+6 \left(2 \cos^2 u  -1\right) \sqrt{2}=0 \\
 27 \left(\frac{2 \left(-2 \cos^2 u  +1 \right) \sqrt{2}}{3}+\sin  u \left(\cos^4 u  -\cos^2 u  -\frac{1}{9}\right) \sqrt{17-9 \cos  4u}\right) \cos  u \sin  u =0
 \end{cases}
 \end{equation}
 \normalsize
 The \textbf{solve} command yields a heavy answer, involving the place holder \textit{RootOf}. Resolving it with the \textbf{allvalues} command is not enough, and the \textbf{evalf} command has to be applied to obtain numerical answers. The output is composed of 2 real values and 6 non real values, which are irrelevant to our purpose. The obtained values are $u=- 2.670495246, 1.900581007,- 0.4710974309, 1.241011647$, which correspond to the 4 points which have been conjectured from Figure \ref{fig envelopes}. The final plot is displayed in Figure \ref{fig envelope with singular points}.
 
 \section{Offsets}
\label{section offsets}
\subsection{The offset as a geometric locus}
\begin{figure}[htb]
\begin{center}
($d=1/3$) \includegraphics[width=4cm]{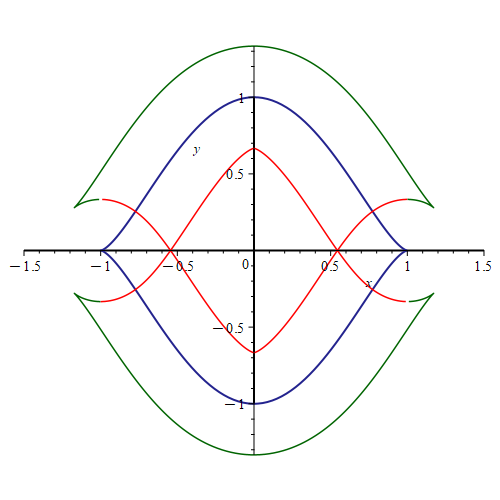}
\qquad
($d=1/2$)\includegraphics[width=4cm]{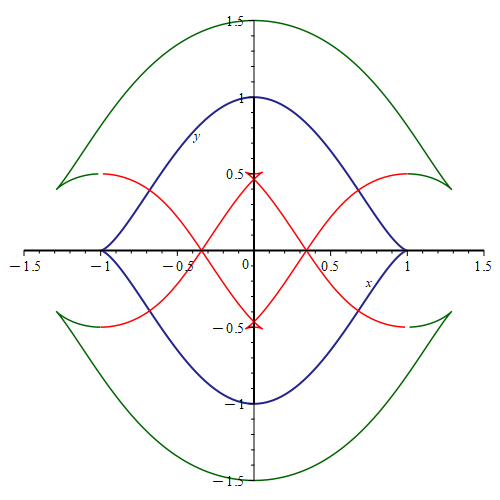}
\qquad
($d=1$) \includegraphics[width=4cm]{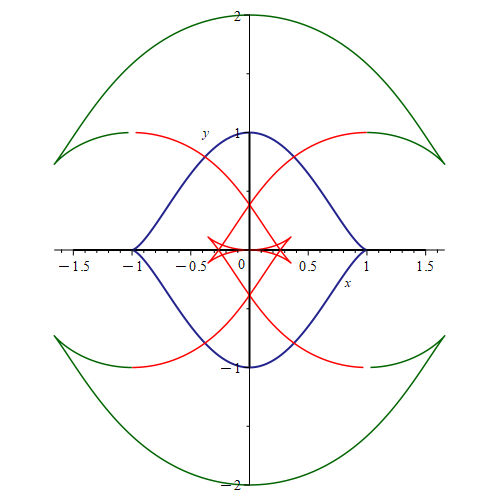}
\qquad
($d=2$)\includegraphics[width=4cm]{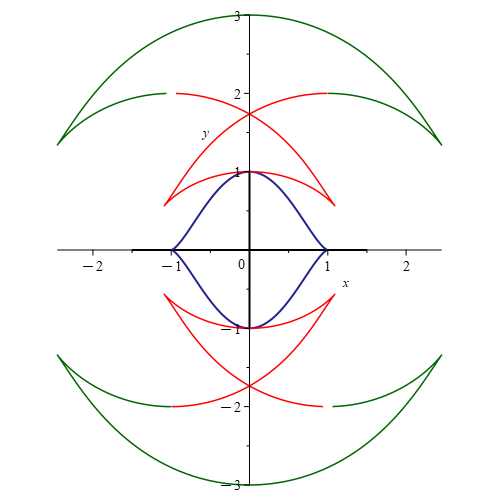}
\caption{Offsets of the kiss curve (the kiss curve is in blue)}
\label{plots offsets}
\end{center}
\end{figure}

Here too, we consider the kiss curve $\mathcal{C}$ given by Equations (\ref{eq kiss curve}). At regular points, a tangent vector is given  by
\begin{equation}
\label{tgt vector}
\overset{\rightarrow}{V}=\left( \frac{dx}{du}, \frac{dy}{du} \right) = \left( -\sin u, 3 \sin^2 u \; \cos u \right). 
\end{equation}
We have here a confirmation that there are singular points for $u=k \pi, \; k \in \mathbb{Z}$. For other values of the parameter, a normal vector is
\begin{equation}
\label{normal vector}
\overset{\rightarrow}{N}= \left( 3 \sin^2 u \; \cos u, \sin u \right). 
\end{equation}
Now we have 
\begin{equation}
\label{norm of tgt and normal vectors}
| \overset{\rightarrow}{V} |= | \overset{\rightarrow}{N} |=  \sqrt{ \sin  \left(u\right)^{2}+9 \sin  \left(u\right)^{4} \cos  \left(u\right)^{2}}.
\end{equation}
and a unit normal vector is 
\begin{equation}
\label{unit normal vector}
\overset{\rightarrow}{N_1}=\frac{1}{\sqrt{\sin  \left(u\right)^{2}+9 \sin  \left(u\right)^{4} \cos  \left(u\right)^{2}}}\; \left( 3 \sin^2 u \; \cos u, \sin u \right).
\end{equation}
Actually, at every regular point there exist 2 unit normal vectors, one pointing inwards and the other one pointing outwards. The internal offset at distance $d$ (with $d>0$) is thus given by the following parametrization
\begin{equation}
\label{param eq external offset}
\begin{cases}
x_I= \cos u + \frac{3d \sin^2 u \cos u}{\sqrt{ \sin^2 u + 9 \sin^4 u \cos^2 u}} \\
y_I=  \sin^3 u + \frac{d \sin u}{\sqrt{ \sin^2 u + 9 \sin^4 u \cos^2 u}}
\end{cases}
\end{equation}
 and the external offset at distance $d$ is given by
 \begin{equation}
\label{param eq internal offset}
\begin{cases}
x_O= \cos u - \frac{3d \sin^2 u \cos u}{\sqrt{ \sin^2 u + 9 \sin^4 u \cos^2 u}} \\
y_O=  \sin^3 u - \frac{d \sin u}{\sqrt{ \sin^2 u + 9 \sin^4 u \cos^2 u}}
\end{cases}
\end{equation}
It is now possible to plot the offsets for various values of the distance $r$. Figure \ref{plots offsets} shows the offsets at distance 1/3, 1/2, 1 and 2.
For $d=2$, it seems that the offset has an arc which coalesces with the kiss curve. This is not true, a suitable zooming will show the difference.

The Maple code for this subsection is as follows:
\small
\begin{verbatim}
setoptions(scaling = constrained):with(plots):
x := cos(u):
y := sin(u)^3:
kiss := plot([x, y, u = 0 .. 2*Pi], color = navy, scaling = constrained, 
         thickness = 2):
axes := plot({[0, t, t = -1 .. 1], [t, 0, t = -1.5 .. 1.5]}, color = black, 
         thickness = 2, labels = ['x', 'y']):
pp := display(axes, kiss):
dx := diff(x, u):
dy := diff(y, u):
len := sqrt(dx^2 + dy^2):
nI := Vector([dy/len, -dx/len]): #internal normal unit vector
xI := d*nI[1] + x:
yI := d*nI[2] + y:
nO := Vector([-dy/len, dx/len]): #external normal unit vector
xO := d*nO[1] + x:
yO := d*nO[2] + y:
for d from 1/3 to 2 by 1/6 do
d;
p1 := plot({[xI, yI, u = 0 .. Pi - 0.01], [-xI, -yI, u = 0 .. Pi - 0.01]}, 
        color = "darkgreen"):
p2 := plot({[xO, -yO, u = 0 .. Pi - 0.01], [-xO, yO, u = 0 .. Pi - 0.01]},
        color = "red"):
display(pp, p1, p2):
end do;
\end{verbatim}
The syntax in the two last plot commands has been chosen in order to avoid irrelevant segments in the figure.

\subsection{Comparison of envelopes of circles centred on the kiss curve and of offsets of the same curve}
\begin{figure}[htb]
\begin{center}
\includegraphics[width=4cm]{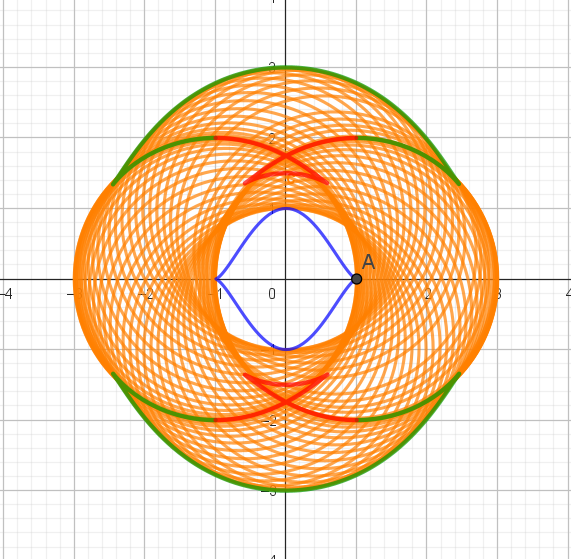}
\qquad
\includegraphics[width=4cm]{envelope_with_singular_points.png}
\caption{Comparison of an envelope and an offset for the same distance}
\label{fig envelope vs offset}
\end{center}
\end{figure}

An applet enabling to compare the envelope and the offset for the same value of $d$ is available at \url{https://www.geogebra.org/m/uzagtrwt}. There, a double animation is available, according to the distance and to the parameter on the curve. Screenshots are displayed in Figure \ref{fig envelope vs offset}: on the left with GeoGebra - the progenitor appears in blue, the 2 components of the offset in red and green, and the shape of the envelope can be understood; on the right with Maple. Note that two circular arcs appear on the left and the right, and 2 other circular arcs in the interior part of the figure. These arcs have a meaning according to another definition of envelopes, describing a larger one than what we have with Definition \ref{def envelope analytic} (see \cite{bruce and giblin}, Chap. 5). Such a situation has been described in \cite{safety}.

\section{Singular points of the offset}
Strong zooming with the GeoGebra applet reveals that the offset has cusps (also called \emph{beaks}, see \cite{alcazar-sendra}), 2nd kind cusps (\emph{thorns} in \cite{alcazar-sendra}) and points of self-intersection (also called \emph{crunodes}). We propose two different approaches for determination of the cusps
\begin{enumerate}
\item Using a parametric representation of the offset, we solve the system of equations $ \frac{dx}{dt}=\frac{dy}{dt}=0$. The corresponding points are the singular points, and we need to classify them.
\item We can also use the curvature to determine the cusps. 
\end{enumerate}
The method to find the crunodes will be different, not requiring differentiation, but also based on the solution of a certain system of equations.

\begin{figure}[htb]
\begin{center}
 \includegraphics[width=6cm]{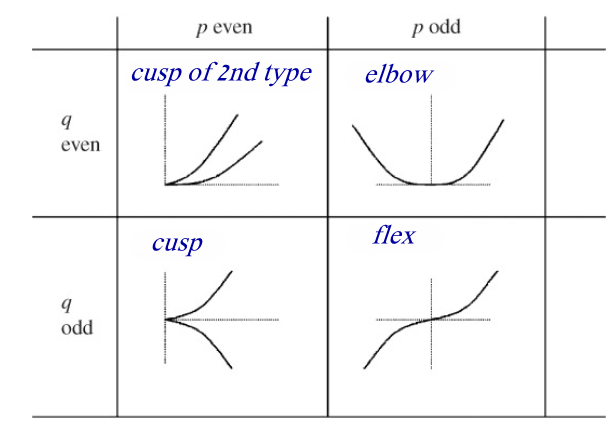}
\caption{Local shapes}
\label{fig local shapes}
\end{center}
\end{figure}

The interactive exploration with GeoGebra leads to the following conjecture:
\begin{itemize}
\item The offset has singular points.
\item At the singular points of the progenitor, no normal vector is defined. This appears as question marks in GeoGebra's algebraic window.
\item The singular points of the offset correspond to regular points of the progenitor (this fact appeared already for the curves studied in \cite{offsets}).
\end{itemize}
 Figure \ref{fig local shapes} classifies the 4 possible shapes.

It remains to check that the obtained values of the parameter correspond to cusps. We refer to \cite{alcazar-sendra} (from whom Figure \ref{fig local shapes} is taken) for a classification of the local shapes of a curve in a neighborhood of a singular point. For any positive integer $n$, let    $\overset{\rightarrow}{V_n}=\left( \frac {d^n x(u)}{du^n},  \frac{d^n y(u)}{du^n} \right)$. Let $p$ be the least value of $n$ such that $\overset{\rightarrow}{V_n}$ is non zero and by $q$ the least value of $n$ such that $\overset{\rightarrow}{V_p}$ and $\overset{\rightarrow}{V_q}$ are linearly independent. Actually, in what follows, we do not check the nature of singular points according to the classification in Figure \ref{fig local shapes}. The reason is that the parametrization of the offset is quite heavy, and the successive derivatives are still heavier (understatement). As the values for which the 1st derivatives vanish are determined only numerically, checking the values of $p$ and $q$ according to the requirements above, is not only unilluminating but will give inaccurate answers. We perform the computations with Maple.

\subsection{Determination of singular points via the vanishing of first derivatives}

As for envelopes, we consider here one component of the offset, WLOG the internal one, given by Equations (\ref{param eq external offset}).
The request to have vanishing 1st derivatives yields the following equations (here for an offset distance equal to 1):
\footnotesize
\begin{equation}
\label{eq 1st derivative offset param}
\begin{cases}
\sin^3 t \; (9 \sqrt{\sin^2 t + 9 \sin^4 t \cos^2 t}\; \cos^4 t + (-9\sqrt{\sin^2 t + 9 \sin^4 t \cos^2 t)} + 6) \; \cos^2 t - \sqrt{\sin^2 t + 9 \sin^4 t \cos^2 t} - 3)=0\\
-3 \cos t \sin^4 t(9\sqrt{\sin^2 t + 9 \sin^4 t \cos^2 t}\; \cos^4 t + (-9 \sqrt{\sin^2 t + 9\sin^4 t\cos^2 t} + 6)\; \cos^2 t - \sqrt{\sin^2 t + 9 \sin^4 t \cos^2 t} - 3)=0
\end{cases}
\end{equation}
\normalsize
This system has  real solutions and  non real (irrelevant) solutions. Numerical values of the real solutions are  2.670495231, 0.4710974228, 1.900581007, 1.241011647. 

\begin{figure}[htb]
\begin{center}
 \includegraphics[width=4cm]{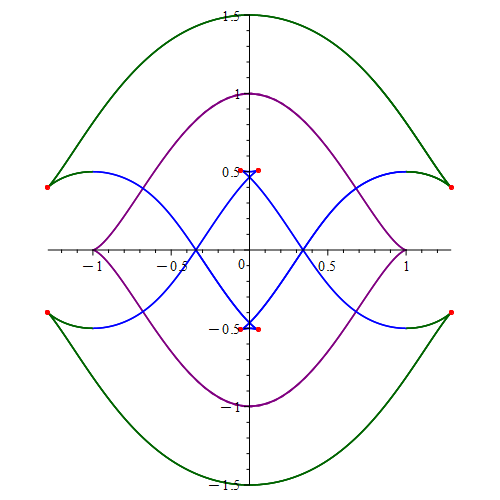}
 \qquad
 \includegraphics[width=4cm]{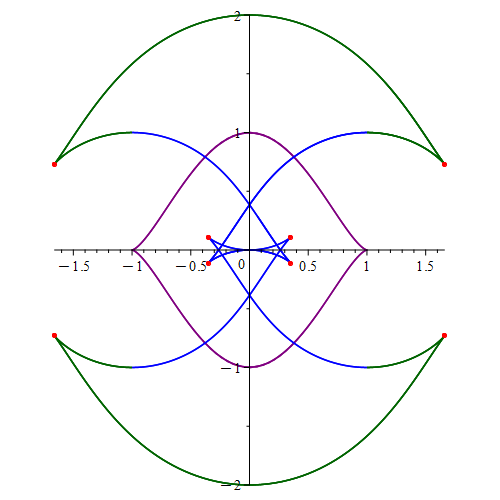}
 \quad
 \includegraphics[width=4cm]{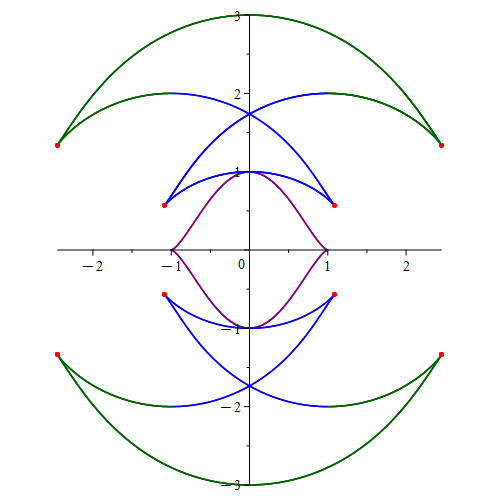} 
\caption{Three examples of offsets with emphasized singular points}
\label{fig offset with emphasized singular points}
\end{center}
\end{figure}

We use the following Maple code:
\small
\begin{verbatim}
x := cos(u);
y := sin(u)^3;
p1 := plot([x, y, u = 0 .. 20], color = "purple");
p2 := plot([x, -y, u = 0 .. 20], color = "purple");
dx := diff(x, u);
dy := diff(y, u);
len := sqrt(dx^2 + dy^2);
nI := Vector([dy/len, -dx/len]);
nO := Vector([-dy/len, dx/len]);
xI := r*nI[1] + x;
yI := r*nI[2] + y;
xO := r*nO[1] + x;
yO := r*nO[2] + y;
r := 0.5;
p3 := plot([xI, yI, u = 0 .. 20], color = "darkgreen",discont=true);
p4 := plot([xO, -yO, u = 0 .. 15], color = "blue",discont=true);
dxi := diff(xI, u);
dxo := diff(xO, u);
sI := solve(dxi = 0, u);
sO := solve(dxo = 0, u);
u_valuesI := [];
for i to nops([sI]) do
    if Im(evalf(sI[i])) = 0 then 
u_valuesI := [op(u_valuesI), evalf(sI[i])]; 
end if;
end do;
u_valuesO := [];
for i to nops([sO]) do
    if Im(evalf(sO[i])) = 0 then 
u_valuesO := [op(u_valuesO), evalf(sO[i])]; end if;
end do;
pointsI_opositeY := [seq([xI(u), yI(u)], u in u_valuesI)];
pointsI_negativeY := [seq([xI(u), -yI(u)], u in u_valuesI)];
pointsO_opositeY := [seq([xO(u), yO(u)], u in u_valuesO)];
pointsO_negativeY := [seq([xO(u), -yO(u)], u in u_valuesO)];
points_values := [op(pointsI_opositeY), op(pointsI_negativeY), 
                 op(pointsO_opositeY), op(pointsO_negativeY)];
points_plot := pointplot(points_values, symbol = solidcircle, symbolsize = 7, 
               color = red);
curve_plotI := plot([xI, yI, u = 0 .. 7], color = blue, thickness = 2,discont=true);
curve_plotO := plot([xO, yO, u = 0 .. 7], color = blue, thickness = 2,discont=true);
curve_plot := plot([x, y, u = 0 .. 2*Pi], color = purple);
display(p1, p2, p3, p4, points_plot, scaling = constrained); 
\end{verbatim}
\normalsize
The solve command yields real and non real answers. The loop which follows filters the answer in order to keep the real ones only\footnote{The authors are aware that this could have been done using a specific option of the solve command, but preferred to do that way in order to have some extra insight, beyond the scope of the present paper.}.  
Figure \ref{fig offset with emphasized singular points} illustrates the answer for respective offset distance $1/2$, 1 and 2: the progenitor is in purple, the offset in green and the cusps highlighted in red.
The crunodes are not emphasized here, we will see that later.

Note the vertical segments in Figure \ref{fig offset with emphasized singular points}; they are superfluous and appear here only because of the way the command \textbf{plot} works. We left them here only to have the unexperienced reader pay attention to the phenomenon. These segments can be cancelled using the \emph{disont=true} option for the \textbf{plot} command, as we will do in later in subsection \ref{subsection self-intersection}. Here the commands should be:
\small
\begin{verbatim}
curve_plotI := plot([xI, yI, u = 0 .. 7], color = blue, thickness = 2,discont=true);
curve_plotO := plot([xO, yO, u = 0 .. 7], color = blue, thickness = 2,discont=true);
\end{verbatim}
\normalsize

\subsection{Determination of cusps using the curvature.}
We use the method described in \cite{curvature}. 
Let $\mathcal{C}$ be a curve determined by a parametric representation $(x(t),y(t))$. If at a point $M(t)=(x(t),y(t))$, both functions $x$ and $y$ are differentiable at least twice, then the curvature of $\mathcal{C}$ at a this point is given by the formula:
\begin{equation}
k(t)=\frac{(r' \times r'') \; e_2}{v^3}=\frac{x'(t)y''(t)-x''(t)y'(t)}{(x'(t)^2+y'(t)^2)^{3/2}},
\end{equation}
where $r(t)=(x(t),y(t))$ and $v=| r'(t) |$.

Define $\hat{k}=\frac{k}{|1+kd|}$; solving the equation: $k=-1/d$ (where the denominator is equal to 0) will give us the $t$ values of the cusps.

\begin{figure}[htb]
\begin{center}
\includegraphics[width=4cm]{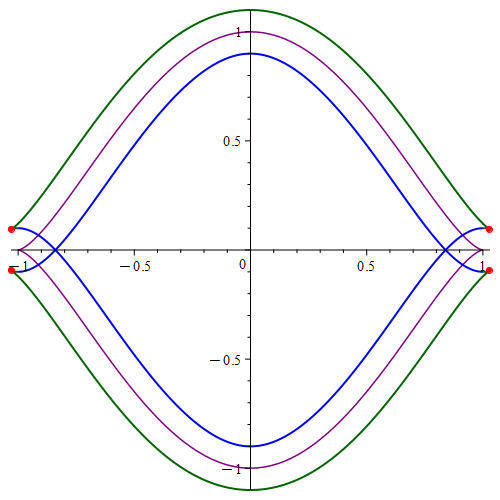}
\qquad
\includegraphics[width=4cm]{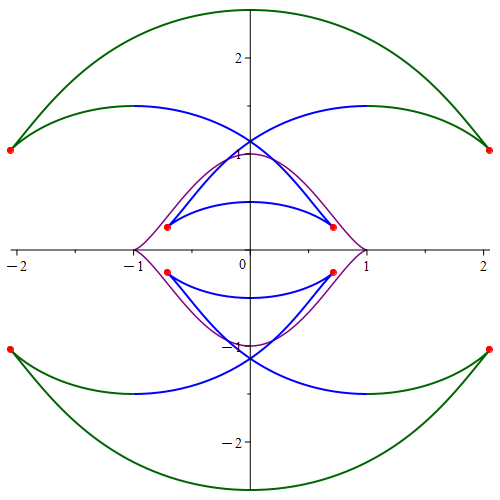}
\qquad
 \includegraphics[width=4cm]{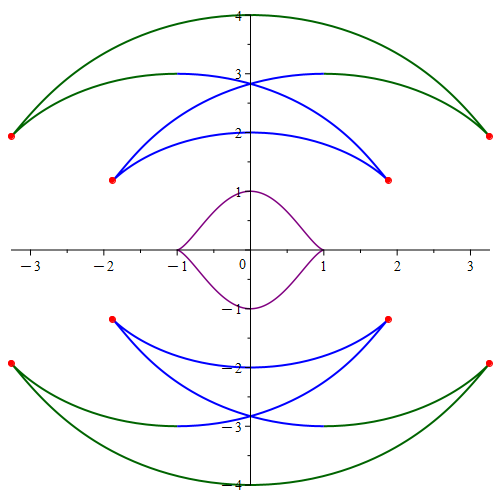}
\caption{Offsets with emphasized singular points (cusps)}
\label{fig offset with emphasized singular points using curvature}
\end{center}
\end{figure}

Here is the implementation in Maple code (note that we are initializing $r$ – for the sake of plotting):
\footnotesize
\begin{verbatim}
x := cos(t);
y := sin(t)^3;
d := 0.1;
dx := diff(x, t);
dy := diff(y, t);
denominator := sqrt(dx^2 + dy^2);
x_offset_pos := x + dy*d/denominator;
y_offset_pos := y - dx*d/denominator;
x_offset_neg := x - dy*d/denominator;
y_offset_neg := y + dx*d/denominator;
r := <x, y>;
dr := diff(r, t);
ddr := diff(dr, t);
kO := (-ddr[1]*dr[2] + ddr[2]*dr[1])/(dr[1]^2 + dr[2]^2)^(3/2);
kI := (ddr[1]*dr[2] - ddr[2]*dr[1])/(dr[1]^2 + dr[2]^2)^(3/2);
sO := solve(kO = -1/d, t);
sI := solve(kI = -1/d, t);
t_values := [];
for i to nops([sO]) do
    if Im(evalf(sO[i])) = 0 then t_values := [op(t_values), evalf(sO[i])]; exp(1)*nd; end if;
end do;
t_values1 := [];
for i to nops([sI]) do
    if Im(evalf(sI[i])) = 0 then t_values1 := [op(t_values1), evalf(sI[i])]; exp(1)*nd; end if;
end do;
outer_singular_points := [seq([evalf(subs(t = t_values[i], x_offset_pos)),
                          evalf(subs(t = t_values[i], y_offset_pos))], i = 1 .. nops(t_values))];
inner_singular_points := [seq([evalf(subs(t = t_values1[i], x_offset_neg)), 
                          evalf(subs(t = t_values1[i], y_offset_neg))], i = 1 .. nops(t_values1))];
curve_plot := plot([x, y, t = 0 .. 2*Pi], color = purple);
outer_singular_points_plot := pointplot(outer_singular_points, symbol = solidcircle, 
                              color = red, symbolsize = 8);
inner_singular_points_plot := pointplot(inner_singular_points, symbol = solidcircle, 
                              color = red, symbolsize = 8);
offset_pos_plot := plot([x_offset_pos, y_offset_pos, t = -Pi .. Pi], 
                              color = darkgreen, thickness = 2, discont = true);
offset_neg_plot := plot([x_offset_neg, y_offset_neg, t = -Pi .. Pi], 
                              color = blue, thickness = 2, discont = true);
display(curve_plot, offset_pos_plot, offset_neg_plot, inner_singular_points_plot, 
                              outer_singular_points_plot);
\end{verbatim}
\normalsize
The result is illustrated in Figure \ref{fig offset with emphasized singular points using curvature} for $d=0.1$, $1.5$ and 3 (from left to right). For $1/2$, 1 and 2, these would have been identical to Figure \ref{fig offset with emphasized singular points}.

We explore these two ways to determine cusps, as an ongoing project addresses a similar question with polynomials of higher degree \cite{Cayley ovals}.
 
\subsection{The points of self-intersection}
\label{subsection self-intersection}
\begin{figure}[htb]
\begin{center}
 \includegraphics[width=4.5cm]{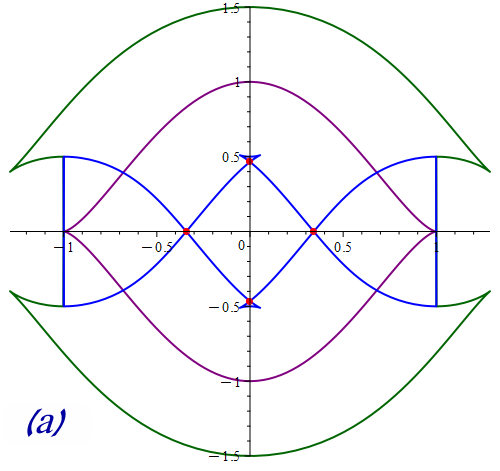}
 \qquad
 \includegraphics[width=4.5cm]{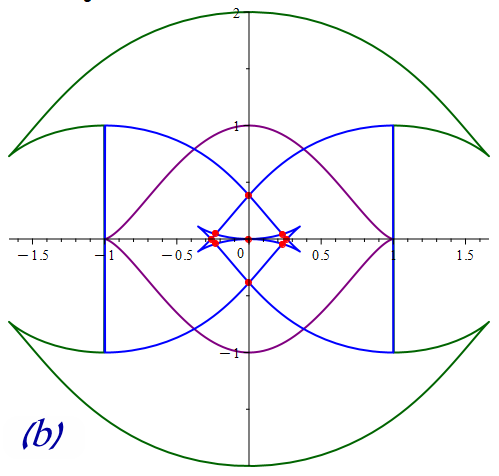}
 \qquad
 \includegraphics[width=4.5cm]{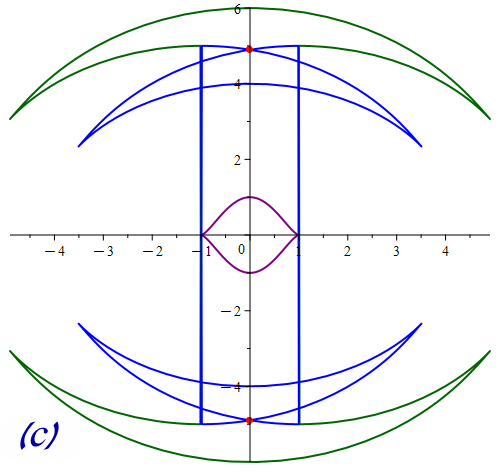}
 \qquad
 \includegraphics[width=4.5cm]{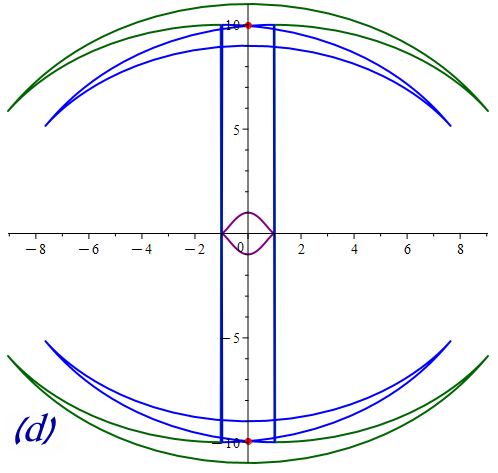}
\caption{An offset with emphasized points of self-intersection on the $x-$axis}
\label{fig crunodes}
\end{center}
\end{figure}

The points of self-intersection correspond to distinct values $s$ and $t$ such that $((x(t),y(t))=(x(s),y(s))$. Actually, we will solve the equations: 
\begin{equation}
\begin{cases}
x(s)+\frac{y'(s) \; d}{\sqrt{x'(s)^2+y'(s)^2}}=x(t)+\frac{y'(t) \; d}{\sqrt{x'(t)^2+y'(t)^2}}\\
y(s)+\frac{x'(s) \; d}{\sqrt{x'(s)^2+y'(s)^2}}=y(t)+\frac{x'(t) \; d}{\sqrt{x'(t)^2+y'(t)^2}}
\end{cases}
\end{equation}
as mentioned in \cite{curvature}. In order to filter the solutions one by one, we need to work in different intervals for the parameter. The determination of the relevant intervals has to be fined tuned in different ways for each offset distance, what reinforces the exploratory aspect of the work. For example, when $d=1$ we have to add a few intervals. Sometimes, an incorrect choice provides irrelevant points.  Finding an automatic way to determine the intervals is still an open issue for the authors.
The Maple code is as follows:
\footnotesize
\begin{verbatim}
x := u -> cos(u);
y := u -> sin(u)^3;
dx := D(x);
dy := D(y);
d := 0.5;
x_offset_pos := u -> x(u) + dy(u)*d/sqrt(dx(u)^2 + dy(u)^2);
y_offset_pos := u -> y(u) - dx(u)*d/sqrt(dx(u)^2 + dy(u)^2);
x_offset_neg := u -> x(u) - dy(u)*d/sqrt(dx(u)^2 + dy(u)^2);
y_offset_neg := u -> y(u) + dx(u)*d/sqrt(dx(u)^2 + dy(u)^2);
eq1 := x_offset_neg(s) - x_offset_neg(t) = 0;
eq2 := y_offset_neg(s) - y_offset_neg(t) = 0;
solutions := [];
ranges := [{s = -4*Pi .. -3*Pi, t = -4*Pi .. -3*Pi}, {s = -(7*Pi)/2 .. -(5*Pi)/2, t = -(7*Pi)/2 .. -(5*Pi)/2}, 
               {s = -(5*Pi)/2 .. -(3*Pi)/2, t = -(5*Pi)/2 .. -(3*Pi)/2}, {s = 3*Pi .. 4*Pi, t = 3*Pi .. 4*Pi}];
for r in ranges do
    sol := fsolve({eq1, eq2}, r);
    if sol <> NULL and evalf(subs(sol, s)) <> evalf(subs(sol, t)) then
        solutions := [op(solutions), sol];
    end if;
end do;
points := [];
for sol in solutions do
    s_sol := evalf(subs(sol, s));
    t_sol := evalf(subs(sol, t));
    x_pos := evalf(x_offset_neg(s_sol));
    y_pos := evalf(y_offset_neg(s_sol));
    points := [op(points), [x_pos, y_pos]];
end do;
plot1 := plot([x(u), y(u), u = -Pi .. Pi], color = purple, thickness = 2);
plot2 := plot([x_offset_pos(u), y_offset_pos(u), u = -Pi .. Pi], 
               color = darkgreen, thickness = 2);
plot3 := plot([x_offset_neg(u), y_offset_neg(u), u = -Pi .. Pi], 
color = blue, thickness = 2);
points_plot := pointplot(points, symbol = solidcircle, color = red, symbolsize = 10);
combined_plot := display([plot1, plot2, plot3, points_plot], 
               title = "Original Curve, Offsets, 
and Solution Points", titlefont = [Helvetica, 14]);
combined_plot;
\end{verbatim}
\normalsize
The output for $d=1/2,1,5,10$ is displayed in Figure \ref{fig crunodes}.

First runs showed irrelevant vertical segments in the various components (a well-known phenomenon); they have been removed using the \emph{discont=true} option in the plot command.

\section{Discussion}

The study of envelopes and offsets is a classic in Differential Geometry. Thom \cite{thom} says that one of the problems of this theory is that it has too few theorems and too many special cases. 
With time, numerous examples have been explored, some of them studied in papers in reference here. Following Thom's claim, it is only natural to explore always new situations. New curves appear, and their topology is an interesting topic. In particular the determination of singular points, cusps, crunodes, etc.

Technology provides important tools for the segment exploration-conjecture-proof. Some kinds of software provide a useful environment for exploration, with interactive features and easy animations. Others may be stronger for the algebraic part, with packages devoted to Gr\"obner bases and Elimination. Original contribution to the study of envelopes is part of the automated methods for Geometry, as in GeoGebra-Discovery, described in \cite{ART} and partly in \cite{DP-ThEdu}. We could not really use this here, as the \textbf{Envelope} command still works for situations more elementary that what we considered in this paper. We wish that in the near future this will be improved.  What we could use for exploration was one of the version of the \textbf{Locus} command and the animations with \emph{Trace On} offered by GeoGebra.

Finally, despite the spirit of the period, we wish to mention that this paper has been written without help from any kind of Generative AI. 

\paragraph{Acknowledgment:} The authors wish to thank Omer Yagel, from DigiSec Ltd. (Israel), for useful suggestions with the Maple code.

\end{document}